\documentclass[12pt, leqno]{amsart}
\usepackage{amssymb,amsfonts,amscd,amsbsy}
\usepackage[mathscr]{eucal}
\theoremstyle{plain}
\newtheorem{thm}{Theorem}[section]
\newtheorem{lem}[thm]{Lemma}

\theoremstyle{definition}

\newtheorem{rem}[thm]{Remark}
\newtheorem{exa}[thm]{Example}

\begin{document}
\title{A note on 4-rank densities}
\author[Robert Osburn]{Robert Osburn}
\address{Department of Mathematics $\&$ Statistics, McMaster University, 1280 Main Street West, Hamilton, Ontario L8S 4K1}
\email{osburnr@icarus.math.mcmaster.ca}
\begin{abstract}
For certain real quadratic number fields, we prove
density results concerning 4-ranks of tame kernels.
We also discuss a relationship between 4-ranks of tame kernels and 4-class ranks of narrow ideal class groups. Additionally, we give a product formula for a
local Hilbert symbol.
\end{abstract}

\maketitle

\section{Introduction}
Let F be a real quadratic number field and $\mathcal{O}_F$ its ring of
integers. In \cite{HK98}, the authors gave an algorithm for
computing the 4-rank of the tame kernel $K_2(\mathcal{O}_F)$. The
idea of the algorithm is to consider matrices with Hilbert symbols
as entries and compute matrix ranks over $\Bbb F_2$. Recently, the author used these matrices to obtain ``density results"
concerning the 4-rank of tame kernels, see \cite{RO}, \cite{MO}.

In this note, we consider the 4-rank of $K_2(\mathcal{O})$ for the real quadratic number fields $\mathbb Q(\sqrt{p_1p_2p_3})$ for primes $p_1 \equiv p_2
\equiv p_3 \equiv 1 \bmod 8$. We will see that

\begin{center}
4-rank $K_2(\mathcal{O}_{\mathbb Q(\sqrt{p_1p_2p_3})}) =$ 0, 1, 2,
or 3.
\end{center}

For squarefree, odd integers d, consider the set

\begin{center}
$X= \{d: d=p_1p_2p_3$, $p_i \equiv 1 \bmod 8 \}$
\end{center}

\noindent for distinct primes $p_i$.

Using GP/PARI \cite{pari}, we computed the following: For $50881
\leq d < 2 \times 10^7$, there are 7257 d's in X. Among them,
there are 2121 d's (29.23$\%$) yielding 4-rank 0, 3977 d's
(54.80$\%$) yielding 4-rank 1, 1086 d's (14.96$\%$) yielding
4-rank 2, and 73 d's (1.01$\%$) yielding 4-rank 3. In fact, we
prove

\begin{thm} \label{T:main}
For the fields $\mathbb Q(\sqrt{p_1p_2p_3})$, 4-rank 0, 1, 2, and
3 appear with natural density $\frac{1}{4}$, $\frac{17}{32}$,
$\frac{13}{64}$, and $\frac{1}{64}$ respectively in $X$.
\end{thm}

In the appendix we point out a beautiful result which may not be
well known. It is a product formula from \cite{San} for a
certain local Hilbert symbol. This product formula both simplifies
numerical computations and is a generalization of Propositions 4.6 and 4.4 in \cite{CH} and \cite{MO}, respectively.

\section{Matrices}

Hurrelbrink and Kolster \cite{HK98} generalize Qin's approach in
\cite{Qin1}, \cite{Qin2} and obtain 4-rank results by computing
$\Bbb F_2$-ranks of certain matrices of local Hilbert symbols.
Specifically, let $F=\mathbb Q(\sqrt{d})$, $d > 1$ and squarefree.
Let $p_1, p_2, \dots , p_t$ denote the odd primes dividing $d$.
Recall 2 is a norm from F if and only if all $p_{i}$'s are $\equiv \pm 1
\bmod 8.$ If so, then d is a norm from $\mathbb Q(\sqrt{2})$, thus
$$
d=u^2-2w^2
$$
for $u,w \in \mathbb Z$. Now consider the matrix: \\

$M_{F/\mathbb Q}=
\left( \begin{matrix} (-d,p_1)_2 & (-d,p_1)_{p_1} \hspace{.1in} \dots & (-d,p_1)_{p_t} \\
(-d,p_2)_2 & (-d,p_2)_{p_1} \hspace{.1in} \dots & (-d,p_2)_{p_t} \\
\vdots & \vdots & \vdots \\
(-d,p_{t-1})_2 & (-d,p_{t-1})_{p_1} \hspace{.1in} \dots & (-d,p_{t-1})_{p_t} \\
(-d,v)_2 & (-d,v)_{p_1} \hspace{.1in} \dots & (-d,v)_{p_t} \\
(d,-1)_2 & (d,-1)_{p_1} \hspace{.1in} \dots & (d,-1)_{p_t} \\
\end{matrix} \right)$. \\

If 2 is not a norm from F, set $v=2$. Otherwise, set $v=u+w$.
Replacing the $1$'s by $0$'s and the $-1$'s by $1$'s, we calculate
the matrix rank over $\mathbb F_2$. From \cite{HK98},

\begin{lem} \label{L:hk} Let $F=\mathbb Q(\sqrt{d})$, $d > 0$ and squarefree. Then
\begin{center}
4-rank $K_2(\mathcal{O}_F) = t$ $-$ rk $(M_{F/\mathbb Q})$ +
$a^{'} - a$
\end{center}

where \\
\begin{center}
$a= \left \{ \begin{array}{l}
0 \quad \mbox{if 2 is a norm from F} \\
1 \quad \mbox{otherwise}
\end{array}
\right. $ \\
\end{center}
and \\

\begin{center}
$a^{'}= \left \{ \begin{array}{l}
0 \quad \mbox{if both -1 and 2 are norms from F} \\
1 \quad \mbox{if exactly one of -1 or 2 is a norm from F} \\
2 \quad \mbox{if none of -1 or 2 are norms from F.} \\
\end{array}
\right. $
\end{center}
\end{lem}

Recall that our case is $\mathbb Q(\sqrt{p_1p_2p_3})$ for primes
$p_1 \equiv p_2 \equiv p_3 \equiv 1 \bmod 8$. In this case
$a=a^{'}$ and we may delete the last row of $M_{F/\mathbb Q}$
without changing its rank (see discussions preceding Proposition
5.13 and Lemma 5.14 in \cite{HK98}). Also note that $v$ is an
$p_1$-adic unit and hence
\begin{center}
$(-p_1p_2p_3, v)_{p_1}=(p_1,v)_{p_1}=\Big(\frac{v}{p_1}\Big)$.
\end{center}
Similarly, $(-p_1p_2p_3,v)_{p_2}= \Big(\frac{v}{p_2}\Big)$ and
$(-p_1p_2p_3,v)_{p_3}= \Big(\frac{v}{p_3}\Big)$. From Lemma 2.1 we
have
\begin{center}
4-rank $K_2(\mathcal{O}_F) = 3$ $-$ rk $(M_{F/\mathbb Q})$
\end{center}
and the matrix $M_{F/\mathbb Q}$ is of the form

\begin{center}
$\left( \begin{matrix} 1 & \Big(\frac{p_2}{p_1}\Big)\Big(\frac{p_3}{p_1}\Big) & \Big(\frac{p_1}{p_2}\Big) & \Big(\frac{p_1}{p_3}\Big)\\
& & & \\
1 & \Big(\frac{p_2}{p_1}\Big) &
\Big(\frac{p_1}{p_2}\Big)\Big(\frac{p_3}{p_2}\Big) &
\Big(\frac{p_2}{p_3}\Big) \\
& & & \\
(-d,u+w)_2 & \Big(\frac{v}{p_1}\Big)  & \Big(\frac{v}{p_2}\Big) &  \Big(\frac{v}{p_3}\Big) \\
\end{matrix} \right)$.
\end{center}

Let us now prove Theorem \ref{T:main}.

\begin{proof}

The idea in \cite{RO} and \cite{MO} is to first consider an
appropriate normal extension $N$ of $\mathbb Q$ and then relate
the splitting of the primes $p_i$ in $N$ to their representation
by certain quadratic forms. The next step is classifying 4-rank
values in terms of values of the symbols $(-d,v)_2$,
$\Big(\frac{v}{p_i}\Big)$. The values of these symbols are then
characterized in terms of $p_i$ satisfying the alluded to
quadratic forms. We then associate Artin symbols to the primes
$p_i$ and apply the Chebotarev density theorem. In what follows,
we classify the 4-rank values in terms of the symbols $(-d,v)_2$,
$\Big(\frac{v}{p_i}\Big)$ and in parenthesis give the relevant
densities in $X$ obtained by using the above machinery. Let us
consider the following four cases (see Table III in \cite{Qin2}).

Case 1: Suppose $\Big(\frac{p_2}{p_1}\Big) = \Big(\frac{p_3}{p_1}\Big) = \Big(\frac{p_3}{p_2}\Big) = 1$. Then we immediately have that
\begin{itemize}
\item 4-rank $K_2(\mathcal{O}_{\mathbb Q(\sqrt{p_1p_2p_3})}) = 3 \iff$ rank $(M_{F/\mathbb Q}) =0 \iff (-d,v)_2=1$ and $\Big(\frac{v}{p_1}\Big)=\Big(\frac{v}{p_2}\Big)=\Big(\frac{v}{p_3}\Big)=1$ ($\frac{1}{64}$) \\

\item 4-rank $K_2(\mathcal{O}_{\mathbb Q(\sqrt{p_1p_2p_3})}) = 2 \iff$ rank $(M_{F/\mathbb Q}) =1 \iff (-d,v)_2=-1$ or $(-d,v)_2=1$ and $\Big(\frac{v}{p_1}\Big)=\Big(\frac{v}{p_2}\Big)=-1$ and $\Big(\frac{v}{p_3}\Big)=1$ or $(-d,v)_2=1$ and $\Big(\frac{v}{p_1}\Big)=\Big(\frac{v}{p_3}\Big)=-1$ and $\Big(\frac{v}{p_2}\Big)=1$ or $(-d,v)_2=1$ and $\Big(\frac{v}{p_2}\Big)=\Big(\frac{v}{p_3}\Big)=-1$ and $\Big(\frac{v}{p_1}\Big)=1$  ($\frac{7}{64}$). \\
\end{itemize}

Case 2: Suppose $\Big(\frac{p_3}{p_2}\Big) = \Big(\frac{p_3}{p_1}\Big) = 1$, $\Big(\frac{p_2}{p_1}\Big) = -1$. Thus
\begin{itemize}
\item 4-rank $K_2(\mathcal{O}_{\mathbb Q(\sqrt{p_1p_2p_3})}) = 2 \iff$ rank $(M_{F/\mathbb Q}) =1 \iff (-d,v)_2=1$ and $\Big(\frac{v}{p_1}\Big)=\Big(\frac{v}{p_2}\Big)=1$ or $(-d,v)_2=1$ and $\Big(\frac{v}{p_1}\Big)=\Big(\frac{v}{p_2}\Big)=-1$ ($\frac{3}{32}$) \\

\item 4-rank $K_2(\mathcal{O}_{\mathbb Q(\sqrt{p_1p_2p_3})}) = 1 \iff$ rank $(M_{F/\mathbb Q}) =2 \iff (-d,v)_2=-1$ or $(-d,v)_2=1$ and $\Big(\frac{v}{p_1}\Big)=\Big(\frac{v}{p_3}\Big)=-1$ and $\Big(\frac{v}{p_2}\Big)=1$ or $(-d,v)_2=1$ and $\Big(\frac{v}{p_2}\Big)=\Big(\frac{v}{p_3}\Big)=-1$ and $\Big(\frac{v}{p_1}\Big)=1$ ($\frac{9}{32}$). \\
\end{itemize}

Case 3: Suppose $\Big(\frac{p_2}{p_1}\Big) = \Big(\frac{p_3}{p_1}\Big) = -1$, $\Big(\frac{p_3}{p_2}\Big) = 1$. Thus
\begin{itemize}
\item 4-rank $K_2(\mathcal{O}_{\mathbb Q(\sqrt{p_1p_2p_3})}) = 1 \iff$ rank $(M_{F/\mathbb Q}) =2 \iff (-d,v)_2=1$ ($\frac{3}{16}$). \\

\item 4-rank $K_2(\mathcal{O}_{\mathbb Q(\sqrt{p_1p_2p_3})}) = 0 \iff$ rank $(M_{F/\mathbb Q}) =3 \iff (-d,v)_2=-1$ ($\frac{3}{16}$) \\
\end{itemize}

Case 4: Suppose $\Big(\frac{p_2}{p_1}\Big) = \Big(\frac{p_3}{p_1}\Big) = \Big(\frac{p_3}{p_2}\Big) = -1$. Then
\begin{itemize}
\item 4-rank $K_2(\mathcal{O}_{\mathbb Q(\sqrt{p_1p_2p_3})}) = 1 \iff$ rank $(M_{F/\mathbb Q}) =2 \iff (-d,v)_2=1$ ($\frac{1}{16}$) \\

\item 4-rank $K_2(\mathcal{O}_{\mathbb Q(\sqrt{p_1p_2p_3})}) = 0 \iff$ rank $(M_{F/\mathbb Q}) =3 \iff (-d,v)_2=-1$ ($\frac{1}{16}$). \\
\end{itemize}

Thus 4-rank 0, 1, 2, and 3 occur with natural density $\frac{1}{16} + \frac{3}{16} = \frac{1}{4}$, $\frac{1}{16} + \frac{3}{16} + \frac{9}{32} = \frac{17}{32}$, $\frac{3}{32} + \frac{7}{64} = \frac{13}{64}$, and $\frac{1}{64}$.

\end{proof}

\begin{rem}
The matrices in \cite{HK98} are related to R\'edei matrices which were used in the 1930's to study the structure of narrow ideal class groups. Namely, for $\mathbb Q(\sqrt{d})$, we considered the case that all odd primes divisors of $d$ are $\equiv 1 \bmod 8$. Thus 2 is a norm from $F=\mathbb Q(\sqrt{d})$ and we have the representation
$$
d=u^2-2w^2.
$$
Let $d^{\prime} = \prod_{i=1}^{t} p_i$. The matrix $M_{F/\mathbb
Q}$ has the form:

\begin{center}
$\left( \begin{matrix} 1 & \\
1 & \\
\vdots & \hat{R}_{F/\mathbb Q} \\
1 & \\
(-d,v)_2 & (-d,v)_{p_1} \hspace{.1in} \dots & (-d,v)_{p_t} \\
\end{matrix} \right)$.
\end{center}

The $(t-1)$ by $t$ matrix $\hat{R}_{F/\mathbb Q}$ can be extended,
without changing its rank, to a $t$ by $t$ matrix $R_{F/\mathbb
Q}$ by adding the last row

\begin{center}
$(-d,p_t)_{p_1}$, $(-d, p_t)_{p_2}, \dots, (-d,p_t)_{p_t}$.
\end{center}

$R_{F/\mathbb Q}$ is known as the R\'edei matrix of the field
$F^{\prime}:= \mathbb Q(\sqrt{d^{\prime}})$ (see \cite{hur} or
\cite{Redei}). Its rank determines the 4-rank of the narrow ideal
class group $C_{F^{\prime}}^{+}$ of the field $F^{\prime}$ by

\begin{center}
4-rank $C_{F^{\prime}}^{+}$ = $t - 1 -$ rank($R_{F/\mathbb Q}$).
\end{center}

Combining this information with Lemma \ref{L:hk}, we have that if
$(-d,u+w)_2=-1$, then 4-rank $K_2(\mathcal{O}_F)$ = 4-rank
$C_{F^{\prime}}^{+}$. Using R\'edei matrices, Gerth in
\cite{Gerth} derived an effective algorithm for computing
densities of 4-class ranks of narrow ideal class groups of
quadratic number fields. It would be interesting to see if
density results concerning 4-class ranks of narrow ideal class
groups (coupled with the product formula in the appendix) can be
used to obtain asymptotic formulas for 4-rank densities of tame
kernels.

\end{rem}

\section{Appendix: A product formula}
Most of the local Hilbert symbols in the martix $M_{F/\mathbb Q}$
are calculated directly. Difficulties arise when $d$ is a norm
from $\mathbb Q(\sqrt{2})$. In this case, we need to calculate the
Hilbert symbols $(-d,u+w)_2$ and $(-d,u+w)_{p_k}$. The local
symbol at 2 is calculated using Lemmas 5.3 and 5.4 in \cite{HK98}.
In this appendix we provide a product formula which allows one to
calculate $(-d,u+w)_{p_k}$ using 2 factors of $d$ at a time.

Let $d$ be a squarefree integer and assume that all odd prime
divisors of $d$ are $\equiv \pm 1 \bmod 8$. Then $d$ is a norm
from $F=\mathbb Q(\sqrt{2})$ and we have the representation
$$
d=u^2-2w^2
$$
with $u>0$. Let $l$ be any odd prime dividing $d$. Note that $l$
does not divide $u+w$ and so

\begin{rem} \label{R:R1}
$(-d, u+w)_l = (l, u+w)_l = \Big(\frac{u+w}{l}\Big)$.
\end{rem}

Recall that any odd prime divisor $l$ of $d$ is $\equiv
\pm 1 \bmod 8$. We fix $x$ and $y$ according to the representation:

\begin{center}
$(-1)^{\frac{l-1}{2}}l = N_{\mathbb Q(\sqrt{2})/\mathbb Q}(x+y\sqrt{2}) = x^2 - 2y^2$
\end{center}

\noindent with $x \equiv 1 \bmod 4$, $x$, $y > 0$. Observe that mod 8, $x$ is odd. Also we can arrange for $x \equiv 1 \bmod 4$ by multiplying $x + y\sqrt{2}$ by $(1+\sqrt{2})^2$.

For $l \equiv 1 \bmod 8$, we have $l=x^2 - 2y^2$ and so $\Big(\frac{l}{y}\Big) =1$. Thus $\Big(\frac{y}{l}\Big) =1$. For $l \equiv 7 \bmod 8$,

\begin{center}
$1= \Big(\frac{-l}{y}\Big) = \Big(\frac{-1}{y}\Big)\Big(\frac{l}{y}\Big) = (-1)^{\frac{y-1}{2}} (-1)^{\frac{y-1}{2}}\Big(\frac{y}{l}\Big) = \Big(\frac{y}{l}\Big)$.
\end{center}

Now let $r$ be an integer not divisible by $l$ which can be
represented as a norm from $\mathbb Q(\sqrt{2})$. Denote by
$\pi_{r} = s + t\sqrt{2}$ an element such that $N_{\mathbb
Q(\sqrt{2})/\mathbb Q}(\pi_{r})=r$ with $s$, $t > 0$. Now let
$u_r$ and $w_r$ be such that

\begin{center}
$u_r + w_r\sqrt{2} = (1+\sqrt{2})(x+y\sqrt{2})(s+t\sqrt{2})$.
\end{center}

By the choice of $x$, $y$, $s$, $t$, we have $u_r > 0$.   Note that

\begin{center}
$N_{\mathbb Q(\sqrt{2})/\mathbb Q}(u_r + w_r\sqrt{2})=-(-1)^{\frac{l-1}{2}}lr$.
\end{center}

Now fix $\frak{l} = <x-y\sqrt{2}>$ a prime ideal above $l$ in
$\mathbb Q(\sqrt{2})$. As $l$ splits in $\mathbb Q(\sqrt{2})$,
$\mathbb Z[\sqrt{2}]/\frak{l} \cong {\mathbb Z}/{l\mathbb Z}$.
This allows us to work mod $\frak{l}$ as opposed to mod $l$. From
the above, $u_r + w_r = 2xs + 3tx + 3sy + 4yt$. Modulo $\frak l$,
we have

\begin{eqnarray}
u_r + w_r & \equiv & 2sy\sqrt{2} + 3ty\sqrt{2} + 3sy + 4yt \nonumber \\
& \equiv & y(3+2\sqrt{2})(s+t\sqrt{2}). \nonumber
\end{eqnarray}

As $\Big(\frac{y}{l}\Big) = 1$, $\Big(\frac{u_r + w_r}{l}\Big) = \Big(\frac{y}{l}\Big)\Big(\frac{{\pi}_r}{\frak{l}}\Big) = \Big(\frac{{\pi}_r}{\frak{l}}\Big)$ where

\begin{center}
$\Big(\frac{{\pi}_r}{\frak{l}}\Big)=
\left \{ \begin{array}{l}
1 \quad \mbox{if $x^2 \equiv {\pi}_r \bmod \frak{l}$ is solvable} \\
-1  \quad \mbox{otherwise}.
\end{array}
\right. $ \\
\end{center}

In the case $r = \prod_{i}^{t-1} p_i$ where $p_i \equiv \pm 1
\bmod 8$, we obtain for each $p_i$ an element $\pi_i \in \mathbb
Q(\sqrt{2})$ of norm $(-1)^{\frac{p_{i} - 1}{2}}p_i$. Let $c$ be
the number of primes dividing $r$ which are congruent to 7 modulo
8. Then we have (up to squares of units) $\pi_r =
(1+\sqrt{2})^{c}\prod_{i}^{t-1}{\pi}_i$. This yields

\begin{center}
$\Big(\frac{u_r + w_r}{l}\Big) = \Big(\frac{(1+\sqrt{2})^c}{\frak{l}}\Big) \prod_{i}^{t-1}\Big(\frac{{\pi}_i}{\frak{l}}\Big)$,
\end{center}

\noindent and so

\begin{center}
$\Big(\frac{u_r + w_r}{l}\Big) = \Big(\frac{u_{-1} + w_{-1}}{l}\Big)^c \prod_{i}^{t-1} \Big(\frac{u_{p_i} + w_{p_i}}{l}\Big)$.
\end{center}

As $-1$ and $2$ are also norms from $\mathbb Q(\sqrt{2})$, we can include $r$'s having factors $-1$ or $\pm 2$. Thus for $r = (-1)^n(2)^m \prod_{i}^{t-1} p_i$ with $m$, $n = 0$, or 1, and each $p_i \equiv \pm 1 \bmod 8$ and $l \neq p_i$ for any $i$, we have

\begin{rem} \label{R:R2}
\begin{center}
$\Big(\frac{u_r + w_r}{l}\Big) = \Big(\frac{u_{-1} + w_{-1}}{l}\Big)^{n+c} \Big(\frac{u_{2} + w_{2}}{l}\Big)^m \prod_{i}^{t-1} \Big(\frac{u_{p_i} + w_{p_i}}{l}\Big)$.
\end{center}
\end{rem}

Setting $r=\frac{d}{l}$, we have $-(-1)^{\frac{l-1}{2}}d = N_{\mathbb Q(\sqrt{2})/\mathbb Q}(u_r + w_r\sqrt{2})$. So for any prime $l \equiv 7 \bmod 8$, we have $N_{\mathbb Q(\sqrt{2})/\mathbb Q}(u_r + w_r\sqrt{2}) = d = N_{\mathbb Q(\sqrt{2})/\mathbb Q}(u + w\sqrt{2})$. Then, up to squares, $\Big(\frac{u_r + w_r}{l}\Big) = \Big(\frac{u + w}{l}\Big)$. For prime divisors $l \equiv 1 \bmod 8$, we have $-d=N_{\mathbb Q(\sqrt{2})/\mathbb Q}(u_r + w_r\sqrt{2})$ and so we include $\Big(\frac{u_{-1} + w_{-1}}{l}\Big)$. To summarize,

\begin{rem} \label{R:R3}
\begin{center}
$(-d,u+w)_l  =
\left \{ \begin{array}{l}
\Big(\frac{u_r + w_r}{l}\Big)\quad \mbox{if $l \equiv 7 \bmod 8$} \\
{\Big(\frac{u_{-1} + w_{-1}}{l}\Big)} \Big(\frac{u_r + w_r}{l}\Big) \quad \mbox{if $l \equiv 1 \bmod 8$}.
\end{array}
\right. $ \\
\end{center}
\end{rem}

We may now reduce to the following $d=rl$: $d=-l$, $d=2l$, and $d=pl$, i.e. calculate the symbols ${\Big(\frac{u_{-1} + w_{-1}}{l}\Big)}$, ${\Big(\frac{u_{2} + w_{2}}{l}\Big)}$, and $\Big(\frac{u_{p} + w_{p}}{l}\Big)$. The first two symbols can be calculated using the following two elementary lemmas.

\begin{lem} \label{L:lem1} ${\Big(\frac{u_{-1} + w_{-1}}{l}\Big)} = 1 \iff (-1)^{\frac{l-1}{2}}l = a^2 - 32b^2$ for some $a$,$b \in \mathbb Z$ with $a \equiv 1 \bmod 4$.
\end{lem}

\begin{lem} \label{L:lem2} ${\Big(\frac{u_{2} + w_{2}}{l}\Big)} = 1 \iff l \equiv \pm 1 \bmod 16$.
\end{lem}

A little care is necessary in computing $\Big(\frac{u_{p} + w_{p}}{l}\Big)$. If ${\Big(\frac{{(-1)}^{\frac{p - 1}{2}}p}{l} \Big)} =1$, then the symbol $\Big(\frac{\pi}{l}\Big)$ is well defined (see discussion preceeding Proposition 3.5 in \cite{CH}) and can be computed using

\begin{lem} For $\mathcal{K}=\mathbb Q(\sqrt{(-1)^{\frac{p - 1}{2}}2p})$ with $p \equiv \pm 1 \bmod 8$ and $h^{+}(\mathcal{K})$ the narrow class number of $\mathcal{K}$, we have

$\Big(\frac{\pi}{l}\Big) = 1 \iff
{l}^{\frac{h^{+}(\mathcal{K})}{4}} = n^{2} - 2pm^{2}$ for some
$n,m \in \mathbb Z$ with $m \not \equiv 0\bmod l$.

For $\mathcal{K}=\mathbb Q(\sqrt{-2p})$ with $p \equiv 7 \bmod 8$,
$\Big(\frac{\pi}{l}\Big) = -1 \iff
{l}^{\frac{h^{+}(\mathcal{K})}{4}} = 2n^{2} + pm^{2}$ for some
$n,m \in \mathbb Z$ with $m \not \equiv 0\bmod l$.

For $\mathcal{K}=\mathbb Q(\sqrt{2p})$ with $p \equiv 1 \bmod 8$, $\Big(\frac{\pi}{l}\Big) = -1 \iff {l}^{\frac{h^{+}(\mathcal{K})}{4}} = {p}n^{2} - 2m^{2}$ for some
$n,m \in \mathbb Z$ with $m \not \equiv 0\bmod l$.
\end{lem}

In fact,

\begin{lem} \label{L:lem3} If ${\Big(\frac{{(-1)}^{\frac{p - 1}{2}}p}{l} \Big)} =1$, then

\begin{center}
$\Big(\frac{u_{p} + w_{p}}{l}\Big) =
\left \{ \begin{array}{l}
\Big(\frac{\pi}{l}\Big)\quad \mbox{for $p \equiv 1 \bmod 8$} \\
{\Big(\frac{u_{-1} + w_{-1}}{l}\Big)} \Big(\frac{\pi}{l}\Big) \quad \mbox{for $p \equiv 7 \bmod 8$}.
\end{array}
\right. $ \\
\end{center}
\end{lem}

The case where ${\Big(\frac{{(-1)}^{\frac{p - 1}{2}}p}{l} \Big)} = -1$ can be done by finding $u_{p}$ and $w_{p}$ from the presentation

\begin{center}
$N_{\mathbb Q(\sqrt{2})/\mathbb Q}(u_{p} + w_{p}\sqrt{2}) = -(-1)^{\frac{p - 1}{2}}{p}{l}$.
\end{center}

Combining Remarks \ref{R:R1}, \ref{R:R2}, and \ref{R:R3}, we have

\begin{thm} \label{T:prod} For $d = (-1)^{n}(2)^{m} \prod_{i=1}^{t} p_i$, with each $p_i \equiv \pm 1 \bmod 8$, we have
\begin{center}
$(-d,u+w)_{p_k} = {\Big(\frac{u_{-1} + w_{-1}}{p_{k}}\Big)}^{n + (-1)^{\frac{p_{k} + 1}{2}}} {\Big(\frac{u_{2} + w_{2}}{p_{k}}\Big)}^m \prod_{i\neq k} \Big(\frac{u_{p_{i}} + w_{p_{i}}}{p_{k}}\Big)$.
\end{center}
\end{thm}

\begin{exa}
Consider the cases $d={\pm}pl, {\pm}2pl$ with $p \equiv 7 \bmod 8$, $l \equiv 1 \bmod 8$, and $\Big(\frac{l}{p}\Big) = 1$ (see Proposition 4.6 in \cite{CH}). Note that $\Big(\frac{\pi}{l}\Big)$ is well defined and so Lemma \ref{L:lem3} is applicable.

For $d=pl$, we have $n=0$, $m=0$ and so
\begin{eqnarray}
(-d,u+w)_l & = & {\Big(\frac{u_{-1} + w_{-1}}{l}\Big)}^{-1} {\Big(\frac{u_{2} + w_{2}}{l}\Big)}^0 {\Big(\frac{u_{-1} + w_{-1}}{l}\Big)} \Big(\frac{\pi}{l}\Big) \nonumber \\
& = & \Big(\frac{\pi}{l}\Big). \nonumber
\end{eqnarray}

For $d=2pl$, we have $n=0$, $m=1$. Thus
\begin{eqnarray}
(-d,u+w)_l & = & {\Big(\frac{u_{-1} + w_{-1}}{l}\Big)}^{-1} {\Big(\frac{u_{2} + w_{2}}{l}\Big)}^1 {\Big(\frac{u_{-1} + w_{-1}}{l}\Big)} \Big(\frac{\pi}{l}\Big) \nonumber \\
& = & \Big(\frac{2+\sqrt{2}}{l}\Big) \Big(\frac{\pi}{l}\Big).
\nonumber
\end{eqnarray}

For $d=-pl$, we have $n=1$, $m=0$. This yields
\begin{eqnarray}
(-d,u+w)_l & = & {\Big(\frac{u_{-1} + w_{-1}}{l}\Big)}^{0} {\Big(\frac{u_{2} + w_{2}}{l}\Big)}^0 {\Big(\frac{u_{-1} + w_{-1}}{l}\Big)} \Big(\frac{\pi}{l}\Big) \nonumber \\
& = & \Big(\frac{1+\sqrt{2}}{l}\Big) \Big(\frac{\pi}{l}\Big).
\nonumber
\end{eqnarray}

Finally, for $d=-2pl$, we have $n=1$, $m=1$. So
\begin{eqnarray}
(-d,u+w)_l & = & {\Big(\frac{u_{-1} + w_{-1}}{l}\Big)}^{0} {\Big(\frac{u_{2} + w_{2}}{l}\Big)}^1 {\Big(\frac{u_{-1} + w_{-1}}{l}\Big)} \Big(\frac{\pi}{l}\Big) \nonumber \\
& = & \Big(\frac{2+\sqrt{2}}{l}\Big)
\Big(\frac{1+\sqrt{2}}{l}\Big) \Big(\frac{\pi}{l}\Big). \nonumber
\end{eqnarray}
\end{exa}

\begin{exa}
Consider the cases $d={\pm}pl$ with $p \equiv l \equiv 1 \bmod 8$,
and $\Big(\frac{l}{p}\Big) = 1$ (see Proposition 4.4 in
\cite{MO}). Again $\Big(\frac{\pi}{l}\Big)$ is well defined and so
Lemma \ref{L:lem3} is applicable.

For $d=pl$, we have $n=0$, $m=0$, and so
\begin{eqnarray}
(-d,u+w)_l & = & {\Big(\frac{u_{-1} + w_{-1}}{l}\Big)}^{-1} {\Big(\frac{u_{2} + w_{2}}{l}\Big)}^0 {\Big(\frac{u_{p} + w_{p}}{l}\Big)} \nonumber \\
& = & \Big(\frac{1+\sqrt{2}}{l}\Big) \Big(\frac{\pi}{l}\Big).
\nonumber
\end{eqnarray}

For $d=-pl$, we have $n=1$, $m=0$. Thus
\begin{eqnarray}
(-d,u+w)_l & = & {\Big(\frac{u_{-1} + w_{-1}}{l}\Big)}^{0} {\Big(\frac{u_{2} + w_{2}}{l}\Big)}^0 {\Big(\frac{u_{p} + w_{p}}{l}\Big)} \nonumber \\
& = & \Big(\frac{\pi}{l}\Big). \nonumber
\end{eqnarray}
\end{exa}

\section*{Acknowledgments}
I would like to thank Jurgen Hurrelbrink for bringing Sanford's
thesis to my attention and Manfred Kolster for many productive
discussions. Additionally, I thank Florence Soriano-Gafiuk at the
Universit\'e de Metz for her hospitality.

\end{document}